\documentclass[11pt]{article}
\usepackage{anysize}
\marginsize{2.0cm}{2.0cm}{2.0 cm}{2.0 cm}
\usepackage{setspace}

\usepackage{latexsym,amsfonts,amsmath}
\usepackage{graphicx}
\usepackage{color}
\usepackage{epstopdf}
\usepackage[titletoc]{appendix}

\def\liminf{\mathop{\underline{\lim}}}
\def\limsup{\mathop{\overline{\lim}}}

\newtheorem{condition}{Condition}[section]{\bfseries}{\itshape}

{\bfseries}{\itshape}

\newtheorem{theorem}{Theorem}[section]{\bfseries}{\itshape}

\newtheorem{corollary}{Corollary}[section]{\bfseries}{\itshape}

\newtheorem{proposition}{Proposition}[section]{\bfseries}{\itshape}

\newtheorem{example}{Example}[section]{\bfseries}{\itshape}

\newtheorem{lemma}{Lemma}[section]{\bfseries}{\itshape}

\newtheorem{remark}{Remark}[section]{\bfseries}{\itshape}

\newtheorem{definition}{Definition}[section]{\bfseries}{\itshape}

{\bfseries}{\itshape}

\begin{document}

\title{A useful technique for piecewise deterministic Markov decision processes}
\author{Xin Guo \thanks{Department of Mathematical Sciences, University of
Liverpool, Liverpool, L69 7ZL, U.K.. E-mail:  X.Guo21@liv.ac.uk.}~ and Yi
Zhang \thanks{Corresponding author. Department of Mathematical Sciences, University of
Liverpool, Liverpool, L69 7ZL, U.K.. E-mail: yi.zhang@liv.ac.uk.}}
\date{}
\maketitle

\par\noindent{\bf Abstract:}
This paper presents with justifications a technique that is useful for the study of piecewise deterministic Markov decision processes (PDMDPs) with general policies and unbounded transition intensities. This technique produces an auxiliary PDMDP from the original one. As to be discussed and claified, the auxiliary PDMDP possesses certain desired properties, which may not be possessed by the original PDMDP. Moreover, the performance measure of any policy in the original PDMDP can be replicated by the auxiliary PDMDP for a large class of performance criteria. As an application, we apply this technique to risk-sensitive PDMDPs with total cost criteria.
\bigskip

 Piecewise deterministic Markov decision processes. Unbounded transition intensities.
\bigskip


\section{Introduction}
This paper concerns the optimal control of piecewise deterministic Markov processes, where the state evolves according to a deterministic and uncontrolled flow between two consecutive jumps, and the transition intensities and post-jump distributions are controlled. Henceforth it will be termed as a piecewise deterministic Markov decision process (PDMDP).

A powerful method of studying PDMDPs is to reduce it to an equivalent discrete-time Markov decision process (DTMDP) by inspecting the PDMDP at each of its jump moments and regarding the (possibly relaxed) control function used during a sojourn time as an action in the DTMDP,  see e.g., \cite{BauerleRieder:2011,Costa:2013,Davis:1993,GuoZhangAMO,Yushkevich:1980,Zhang:2017}. Consequently, the action space in the induced DTMDP, as a set of measurable mappings, is in general a more complicated object than the action space in the original PDMDP. The reason for applying this reduction is to gain access to the rich toolbox of known results on DTMDPs that have been studied since 1950s.

It is well appreciated that the theory of DTMDPs is better established when the underlying DTMDP model satisfies some compactness-continuity conditions, see \cite{Hernandez-Lerma:1996,Hernandez-Lerma:1999,Jaskiewicz:2008}. One example of such compactness-continuity conditions is that the action space is a compact Borel space, the loss function is lower semicontinuous in the action, and the transition kernel possesses a strong Feller property with respect to the action, i.e., it maps each bounded measurable function on the state space to a function, which is jointly measurable in the state and action, and also continuous in the action. (A precise formulation of this is in the appendix.)

However, even if the original PDMDP satisfies a natural set of compactness-continuity conditions, see Condition \ref{GZyExponentialCondition02} below, it can happen that the transition kernel in the induced DTMDP fails to satisfy the desired continuity condition. We demonstrate this in Example \ref{JapanExample01} below. On the other hand, it turns out that this inconvenience does not appear if the transition intensities of the PDMDP are strongly positive, i.e., bounded away from zero by a constant, see Proposition \ref{09June2020PropositionA01}.

In the more specific context of continuous-time Markov decision processes (CTMDPs), where the state does not change between two consecutive jumps, another difficulty associated with nonstationary policies and transition intensities not strongly positive, which is different from the aforementioned one, was documented in \cite{Feinberg:2004,GuoXPZhang:2017,Piunovskiy:2015}.

The contribution of this paper is that we present a technique, which produces an auxiliary PDMDP model satisfying the following: a) the original PDMDP can be thinned out of the auxiliary PDMDP, and the performance of any policy in the original PDMDP is replicated by a corresponding policy in the auxiliary PDMDP; b) the transition intensities in the auxiliary PDMDP are strongly positive, and its induced DTMDP satisfies the desired compactness-continuity conditions\footnote{See Proposition \ref{09June2020PropositionA01} for the precise definition of the compactness-continuity condition for the induced DTMDP.} if so does the original PDMDP. Then as an application, we extend some optimality results for risk-sensitive PDMDPs with total cost criteria, which were previously obtained in \cite{GuoZhangAMO} under the extra requirement on the transition intensities being strongly positive, see the footnote in Proposition \ref{09June2020PropositionA01}. This requirement is omitted here with the help of the proposed technique.

In the context of CTMDPs, there have been known techniques, which produce new CTMDP models, in which the performance measure in the original model can be replicated. One example is the uniformization in \cite{puterman,Serfozo:1979} requiring bounded transition intensities and justified under stationary policies. Its more recent variants that deal with models with unbounded transition intensities can be found in \cite{Cekyay:2018,Piunovskiy:2015}, which are also similar in nature: the former one considers CTMDPs with a denumerable state space and deterministic stationary policies, whereas in the latter one, models with more general policies were considered. Our technique can be viewed as their extension from CTMDPs to PDMDPs, but there are some notable differences.

In greater detail, our technique is closer to \cite{Piunovskiy:2015} as we also consider general policies and unbounded transition intensities. Nevertheless, apart from a more general process under control being dealt with here, let us mention another more important difference from \cite{Piunovskiy:2015} as follows. The justification in \cite{Piunovskiy:2015} only serves the risk-neutral problem with total cost criteria, as it is based on the comparison of (total) occupation measures in the original model and the new model. This is not suitable for risk-sensitive problems, where in general the performance measures cannot be readily written as integrals of the cost rate with respect to the occupation measures. In contrast, here our justification applies to both risk-neutral and risk-sensitive problems with both total cost and average cost criteria, see Remark \ref{09June2020Remark01}(c).

The rest of this paper is organized as follows. In Section \ref{NoteGZSec01} we describe the PDMDP model. In Section \ref{NoteGZSec01} we present and prove the main result with an application to risk-sensitive PDMDPs. The appendix presents the induced DTMDP. While the reference to the appendix can be avoided in the main text, it clarifies and demonstrates the issues mentioned in the beginning of this introduction.

\section{Description of PDMDP model}\label{NoteGZSec01}


Let $S$ be a nonempty Borel state space endowed with the Borel $\sigma$-algebra ${\cal B}(S)$, $A$ be
a nonempty Borel action space, and $q$
stand for a signed kernel $q(dy|x,a)$ on ${\cal{B}}(S)$ given
$(x,a)\in S\times A$ such that
\begin{eqnarray*}
\tilde{q}(\Gamma_S|x,a):=q(\Gamma_S\setminus\{x\}|x,a)\ge 0
\end{eqnarray*}
for all $\Gamma_S\in{\cal{B}}(S)$, and
\begin{eqnarray*}
q(S|x,a)=0,~\bar{q}_x=\sup_{a\in A}q_x(a)<\infty,
\end{eqnarray*}
where $q_x(a):=-q(\{x\}|x,a)$ is the transition intensity. The signed kernel $q$ is also called the
transition rate. Between two consecutive jumps, the state of the process evolves according to a measurable mapping $\phi$ from $S\times [0,\infty)$ to $S$, see (\ref{GGGAMO01}) below. It is assumed that for each $x\in S$
\begin{eqnarray}\label{09June2020Eqn08}
\phi(x,t+s)=\phi(\phi(x,t),s),~\forall~s,t\ge 0;~\phi(x,0)=x,
\end{eqnarray}
and $t\rightarrow \phi(x,t)$ is continuous.  Unless stated otherwise, the term of
measurability is always understood in the Borel sense.
Finally let the cost rate $c$ be a $[0,\infty)$-valued measurable function on $S\times A$.

For the rest of this paper, it is convenient to introduce the
following notations. Let $\mathbb{P}(A)$ be the space of
probability measures on ${\cal B}(A)$, endowed with the
standard weak topology. For each $\mu\in \mathbb{P}(A)$,
\begin{eqnarray*}
&& q_x(\mu):=\int_A q_x(a)\mu(da),~\tilde{q}(dy|x,\mu):= \int_A \tilde{q}(dy|x,a)\mu(da),~c(x,\mu):=\int_A c(x,a)\mu(da).
\end{eqnarray*}

\begin{condition}\label{GGZyExponentialConditionExtra}
For each $x\in S$, $\int_{0}^t\overline{q}_{\phi(x,s)}ds<\infty$, and $\int_{0}^t \sup_{a\in A} c(\phi(x,s),a)ds<\infty$, for each $t\in[0,\infty).$
\end{condition}
Condition \ref{GGZyExponentialConditionExtra} is assumed to hold throughout this paper. The integrals in Condition \ref{GGZyExponentialConditionExtra} are well defined because the integrands are universally measurable.

Now we briefly describe the PDMDP with the system primitives $\{S,A,q,\phi,c\}$. Let us take the sample space $\Omega$ by adjoining to the
countable product space $S\times((0,\infty)\times S)^\infty$ the
sequences of the form
$(x_0,\theta_1,\dots,\theta_n,x_n,\infty,x_\infty,\infty,x_\infty,\dots),$
where $x_0,x_1,\dots,x_n$ belong to $S$,
$\theta_1,\dots,\theta_n$ belong to $(0,\infty),$ and
$x_{\infty}\notin S$ is the isolated point. We equip $\Omega$ with
its Borel $\sigma$-algebra $\cal F$.

Let $t_0:=0=:\theta_0,$ and for each $n\geq 0$, and each
element $\omega:=(x_0,\theta_1,x_1,\theta_2,\dots)\in \Omega$, let \begin{eqnarray*}
h_n:=(x_0,\theta_1,\dots,\theta_n,x_n);~
t_n:=t_{n-1}+\theta_n;~ t_\infty(\omega):=\lim_{n\rightarrow\infty}t_n.
\end{eqnarray*}
Then, $(\Omega,{\cal F})$ is the canonical sample space of the marked point process $(t_n,x_n)$ with the mark space $S$, and $\theta_{n}=t_n-t_{n-1}$ is the sojourn time, where the convention of $\infty-\infty:=\infty$ is in use.
Define the process, which evolves according to the flow $\phi$ during a sojourn time: \begin{eqnarray}\label{GGGAMO01}
\xi_t=\left\{\begin{array}{ll} \phi(x_n,t-t_n), & \mbox{ if } t_n\le t<t_{n+1};\\
x_\infty, & \mbox{ if } t_\infty\le
t, \end{array}\right.
\end{eqnarray}
where $x_\infty\notin S$ is an isolated cemetery point. The process is controlled through its local characteristics as follows.

A policy $\pi$ is given by a
sequence $(\pi_n)$ such that, for each $n=0,1,2,\dots,$
$\pi_n(da|h_n,s)$ is a stochastic kernel on
$A$ given $h_n,s$ with $s>0$, and for each $\omega=(x_0,\theta_1,x_1,\theta_2,\dots)\in
\Omega$, $t> 0,$
\begin{eqnarray}\label{09June2020Eqn12}
\pi(da|\omega,t)&:=&I\{t\ge t_\infty\}\delta_{a_\infty}(da)+
\sum_{n=0}^\infty I\{t_n< t\le
t_{n+1}\}\pi_{n}(da|h_n, t-t_n),
\end{eqnarray}
defines a $\mathbb{P}(A\cup\{a_\infty\})$-valued (relaxed) control process,
where $a_\infty\notin A$ is some isolated point. If for some measurable mapping $\varphi$ from $S$ to $A$ such that $\pi_n(da|x_0,\theta_1,\dots,\theta_n,x,t)\equiv\delta_{\varphi(x)}(da)$, then the policy $\pi=(\pi_n)$ is called deterministic stationary and is identified with the mapping $\varphi.$

A policy $\pi$ and an initial state $x$ define a probability measure $P_x^\gamma$ on the canonical sample space, under which $P_x^\pi(x_0=x)=1$, and the conditional distribution of $(\theta_{n+1},x_{n+1})$ given $h_n$ satisfies
\begin{eqnarray}\label{09June2020Eqn01}
&&P_\gamma^\pi(\theta_{n+1}\in \Gamma_1,~x_{n+1}\in \Gamma_2|h_n)\nonumber\\
&=&\int_{\Gamma_1}e^{-\int_0^t \int_A q_{\phi(x_n,s)}(a)\pi_n(da|h_n,s)ds}\int_{A}\tilde{q}(\Gamma_2|\phi(x_n,t),a)\pi_n(da|h_n,t)dt,\nonumber\\
&&~\forall~\Gamma_1\in{\cal B}((0,\infty)),~\Gamma_2\in{\cal B}(S); \\
&&P_\gamma^\pi(\theta_{n+1}=\infty,~x_{n+1}=x_\infty|h_n)=e^{-\int_0^\infty  \int_A q_{\phi(x_n,s)}(a)\pi_n(da|h_n,s)ds}\nonumber
\end{eqnarray}
on $\{x_n\in S\}$.

The proposed technique in this paper will be applied to the risk-sensitive optimal control problem for the PDMDP with a total cost criterion, which reads
\begin{eqnarray*}
&\mbox{Minimize over all policies $\pi$:~}& E_x^\pi\left[e^{\int_0^\infty \int_A
c(\xi_t,a)\pi(da|\omega,t)dt}\right]\\
&&=E_x^\pi\left[e^{\sum_{n=0}^\infty \int_{t_n}^{t_{n+1}} \int_A
c(\phi(x_n,s),a)\pi_{n}(da|h_n,s-t_{n})ds}\right]=:V(x,\pi)
\end{eqnarray*}
Here $\int_{t_n}^{t_{n+1}}$ is understood as $\int_{(t_n,t_n]\bigcap \mathbb{R}}$, and we put $c(x_\infty,a)\equiv 0$.
The value function is defined by $V^\ast(x)=\inf_{\pi}V(x,\pi)$ for all $x\in S.$
We shall call the above system primitives $\{S,A,q,\phi,c\}$ and the corresponding optimal control problem the ``original model'', to distinguish it from the auxiliary model that will appear later.

However, our technique is also applicable to PDMDPs with other performance measure. For instance, one may consider the expected long run average cost defined by
\begin{eqnarray}\label{09June2020Eqn15}
&&\overline{V}(x,\pi):=\limsup_{T\rightarrow \infty} E_x^\pi\left[ \frac{\int_0^T  \int_A
c(\xi_t,a)\pi(da|\omega,t)dt }{T}\right]\nonumber\\
&=&\limsup_{T\rightarrow \infty}  E_x^\pi\left[\frac{\sum_{n=0}^\infty \int_{t_n\wedge T}^{t_{n+1}\wedge T}  \int_A
c(\phi(x_n,s),a)\pi_{n}(da|h_n,s-t_{n})ds }{T}\right],
\end{eqnarray}
where $t_n\wedge T:=\min\{t_n,T\}$. See Remark \ref{09June2020Remark01}(c).

\section{Main result}

Fix $\lambda>0$ in what follows. We introduce an auxiliary model  $\{\breve{S},A,\breve{q},\breve{\phi},\breve{c}\}$ defined in terms of the system primitives of the original model. When there is a danger of confusion, we shall primarily use breves to signify the auxiliary model. Without special explanations, all the objects signified with breves are understood similarly to their counterparts without breves.

Roughly speaking, the auxiliary model arises from inserting additional inspections of the state process during each sojourn time in the original model (up to the moment of explosion) taking place in an independent Poisson process with rate $\lambda$. The changes in the second coordinate of the state in the auxiliary model take place at and only at each of such inspection epochs, which will be recorded as ``fictitious'' jumps and generate strongly positive transition intensities.

The state space is $\breve{S}=S\times \{-1,1\}$, endowed with the product topology, where $\{-1,1\}$ is with the discrete topology. The action space is $A$. The transition rate $\breve{q}$ on ${\cal B}(\breve{S})$ given $\breve{S}\times A$ is defined as follows:
\begin{eqnarray*}
\breve{q}(dy\times\{-i\}|(x,i),a)=\lambda \delta_x(dy);~\breve{q}(dy\times\{i\}|(x,i),a)=q(dy|x,a)-\lambda\delta_x(dy),
\end{eqnarray*}
with $\delta_{x}(dy)$
being the Dirac measure concentrated on the singleton $\{x\},$
so that
\begin{eqnarray*}
&&\breve{q}_{(x,i)}(a)=\breve{q}(\breve{S}\setminus \{(x,i)\}|(x,i),a)=\breve{q}(S\setminus\{x\}\times\{-1,1\}|(x,i),a)+\breve{q}(S\times\{-i\}|(x,i),a)\\
&=&q_x(a)+\lambda~\forall~(x,i)\in \breve{S},~a\in A.
\end{eqnarray*}
In other words, the auxiliary model has strongly positive transition intensities. The flow is defined by \begin{eqnarray*}
\breve{\phi}((x,i),t)=(\phi(x,t),i).
\end{eqnarray*} The cost rate is \begin{eqnarray*}
\breve{c}((x,i),a)=c(x,a)~\forall~(x,i)\in \breve{S},~a\in A.
\end{eqnarray*}

Let
 \begin{eqnarray*}
 \breve{V}((x,i),\breve{\pi})=\breve{E}^{\breve{\pi}}_{(x,i)}[e^{\int_0^\infty \breve{c}(\breve{\xi}_t,a)\breve{\pi}(da|\breve{\omega},t)dt}].
 \end{eqnarray*}

\begin{definition}\label{9thJuneDefinition01}
Consider the canonical sample space of the marked point process $(\breve{t}_n,x_n,i_n)$, and a sample path
\begin{eqnarray*}
\breve{\omega}=((x_0,i_0),\breve{\theta}_1,(x_1,i_1),\breve{\theta}_2,\dots, (x_{n-1},i_{n-1}),\breve{\theta}_n,(x_n,i_n),\dots).
\end{eqnarray*}
We say a mark $(x_l,i_l)$ ($l\ge 1$) is immediately after a fictitious jump if $i_l=-i_{l-1}$, or equivalently, $x_l=\phi(x_{l-1},\breve{\theta}_l)$, where $\breve{\theta}_l$ is the sojourn time before the mark $(x_l,i_l)$. A mark that is not immediately after a fictitious jump is called immediately after an honest jump. We regard $(x_0,i_0)$ as a mark immediately after an honest jump.
\end{definition}
Using the notation in the above definition, we may consider out of $(\breve{t}_n,x_n,i_n)$ another marked point process $(\tau_{(m)},x_{(m)},i_{(m)})$ with $\tau_{(0)}:=0$ by counting only the points with marks immediately after honest jumps. Since $(x_0,i_0)$ is regarded as a mark immediately after an honest jump, $x_0=x_{(0)}$ and $i_0=i_{(0)}$. Since $i_{(m)}=i_{(0)}$ for all $m\ge 0$ almost surely in $(\tau_{(m)},x_{(m)},i_{(m)})$, with $i_{(0)}$ being fixed we may simply consider the marked point process $(\tau_{(m)},x_{(m)})$ instead of $(\tau_{(m)},x_{(m)},i_{(m)})$.

\begin{theorem}\label{09June2020Thm09}
Suppose Condition \ref{GGZyExponentialConditionExtra} is satisfied. For each policy $\pi=(\pi_n)$ in the original PDMDP model, there is a policy $\breve{\pi}=(\breve{\pi}_n)$ in the auxiliary PDMDP model such that for all $x\in S$ and $i\in\{-1,1\}$:
\begin{itemize}
\item[(a)]
The distribution of the marked point process  $(\tau_{(m)},x_{(m)})$ under $\breve{P}_{(i,x)}^{\breve{\pi}}$ coincides with the distribution of the marked point process $(t_m,x_m)$ under $P_{x}^\pi$. In other words, the marked point process in the original model (under $P_x^\pi$) may be thinned out of $(\breve{t}_n,x_n,i_n)$ in the auxiliary model (under $\breve{P}_{(x,i_0)}^{\breve{\pi}}$) by counting only the points with marks immediately after honest jumps.
\item[(b)]
 $V(x,\pi)=\breve{V}((x,i),\tilde{\pi})$.
\end{itemize}
\end{theorem}

\par\noindent\textit{Proof.} We will make use of the notation in Definition \ref{9thJuneDefinition01} freely.

(a) Let a policy $\pi=(\pi_n)$ for the original model be fixed. Consider the corresponding policy $\breve{\pi}=(\breve{\pi}_n)$ in the auxiliary model defined as follows.  For the $n$-history \begin{eqnarray*}
\breve{h}_n=((x_0,i_0),\breve{\theta}_1,(x_1,i_1),\breve{\theta}_2,\dots, (x_{n-1},i_{n-1}),\breve{\theta}_n,(x_n,i_n))
\end{eqnarray*}
in the auxiliary model, let $m=m(\breve{h}_n)$ be the number of honest jumps over $(0,t_n]$ within $\breve{h}_n$, so that if we count the initial mark $(i_0,x_0)$ as immediately after an honest jump, then there are $m+1$ marks immediately after honest jumps within $\breve{h}_n$. Then we define
\begin{eqnarray}\label{09June2020Eqn10}
&&\breve{\pi}_n(da|\breve{h}_n,t)=\pi_m(da|x_0,\tau_{(1)},x_{(1)},\tau_{(2)}-\tau_{(1)},\dots,\tau_{(m)}-\tau_{(m-1)},x_{(m)},t+\breve{t}_n-\tau_{(m)})
\end{eqnarray}
for $t>0.$

Consequently, for each $n,m\ge 0$ and for each $t\in (0,\infty)$ satisfying $t\in (\breve{t}_n,\breve{t}_{n+1}] \subseteq (\tau_{(m)},\tau_{(m+1)}]$, we have
\begin{eqnarray}\label{Correction05}
\breve{\pi}(da|\breve{\omega},t)&=&\breve{\pi}_n(da|\breve{h}_n,t-\breve{t}_n)\nonumber\\
&=&\pi_m(da|x_0,\tau_{(1)},x_{(1)},\tau_{(2)}-\tau_{(1)},\dots,\tau_{(m)}-\tau_{(m-1)},x_{(m)},t-\tau_{(m)}),
\end{eqnarray}
where the first equality is by (\ref{09June2020Eqn12}) applied to $\breve{\pi}$.

For brevity, below we put
\begin{eqnarray}\label{09June2020Eqn05}
&&\tilde{q}(dy|\phi(x_{(m)},t),\pi_m):=
\int_A\tilde{q}(dy|\phi(x_{(m)},t),a)\pi_m(da|x_{(0)},\tau_{(1)},x_{(1)},\tau_{(2)}-\tau_{(1)},\dots,x_{(m)},t); \nonumber\\
&& {q}(dy|\phi(x_{(m)},t),\pi_m):=
\int_A {q}(dy|\phi(x_{(m)},t),a)\pi_m(da|x_{(0)},\tau_{(1)},x_{(1)},\tau_{(2)}-\tau_{(1)},\dots,x_{(m)},t);\nonumber\\
&& q_{\phi(x_{(m)},t)}(\pi_m):=\tilde{q}(S|\phi(x_{(m)},t), \pi_m).
\end{eqnarray}

Now let us show that the distribution of the marked point process $(\tau_m,x_m)$ under $\breve{P}_{(x,i_0)}^{\breve{\pi}}$ coincides with the distribution of the marked point process in the original model under $P_x^\pi$.
To this end, in view of (\ref{09June2020Eqn01}), $x_{(0)}=x_0$ and $\tau_{(0)}=0$, it is sufficient to show that
\begin{eqnarray}\label{09June2020Eqn02}
&&\breve{P}_{(x,i)}^{\breve{\pi}}(x_{(m+1)}\in \Gamma,~\tau_{(m+1)}-\tau_{(m)}\in[0,T]|x_{(0)},\tau_{(1)},x_{(1)},\dots,\tau_{(m)}-\tau_{(m-1)},x_{(m)})\nonumber\\
&=& \int_{0}^T  \tilde{q}(\Gamma|\phi(x_{(m)},t),\pi_m)  e^{-\int_0^t q_{\phi(x_{(m)},s)}(\pi_m)ds}dt
\end{eqnarray}
on $\{\tau_{(m)}<\infty\}$ for each $T>0,$ $\Gamma\in{\cal B}(S)$ and $m\ge 0.$
Note that
\begin{eqnarray*}
&&\breve{P}^{\breve{\pi}}_{(x,i)}(x_{(m+1)}\in \Gamma,~\tau_{(m+1)}-\tau_{(m)}\in[0,T]|~x_{(0)},\tau_{(1)},x_{(1)},\dots,\tau_{(m)}-\tau_{(m-1)},x_{(m)})\\
&=&\sum_{n=0}^\infty \breve{P}^{\breve{\pi}}_{(x,i)}(x_{(m+1)}\in \Gamma,~\tau_{(m+1)}-\tau_{(m)}\in[0,T],\\
&&~ \mbox{exactly $n$ ficticious jumps over $[\tau_{(m)},\tau_{(m+1)}]$}|~x_{(0)},\tau_{(1)},x_{(1)},\dots,\tau_{(m)}-\tau_{(m-1)},x_{(m)}).
\end{eqnarray*}
Since $\sum_{n\ge 0}\frac{\lambda^nt^n}{n!}=e^{\lambda t}$, equality (\ref{09June2020Eqn02}) would be justified once we show that
\begin{eqnarray}\label{09June2020Eqn07}
&&\breve{P}^{\breve{\pi}}_{(x,i)}(x_{(m+1)}\in \Gamma,~\tau_{(m+1)}-\tau_{(m)}\in[0,T],\nonumber\\
&&~ \mbox{exactly $n$ ficticious jumps over $[\tau_{(m)},\tau_{(m+1)}]$}|~x_{(0)},\tau_{(1)},x_{(1)},\dots,\tau_{(m)}-\tau_{(m-1)},x_{(m)})\nonumber\\
&=& \int_0^T \frac{\lambda^n t^n}{n!}\tilde{q}(\Gamma|\phi(x_{(m)},t),\pi_m)e^{-\lambda t}e^{-\int_0^t q_{\phi(x_m,s)}(\pi_m) ds}.
\end{eqnarray}
For this, let us verify for each $T>0$ and $y\in S$ that
\begin{eqnarray}\label{09June2020Eqn06}
&&f(T,y,n-1):=\int_0^T \int_0^{T-r_1}\int_{0}^{T-\sum_{i=1}^2r_i}\dots\int_{0}^{T-\sum_{i=1}^{n-1}r_i} \tilde{q}(\Gamma|\phi(y,\sum_{i=1}^{n-1}r_i+t),\pi_m)\nonumber\\
&&\times e^{-\int_0^t (q_{\phi(y,\sum_{i=1}^{n-1}r_i+s)}(\pi_m)+\lambda)ds}\lambda^{n-1} e^{-\int_0^{\sum_{i=1}^{n-1}r_i} (q_{\phi(y,s)}(\pi_m)+\lambda)ds}   dt dr_{n-1}\dots d r_1\nonumber\\
&=& \int_0^T \frac{\lambda^{n-1}t^{n-1}}{(n-1)!} \tilde{q}(\Gamma|\phi(y,t),\pi_m)e^{-\int_{0}^t (q_{\phi(y,s)}(\pi_m)+\lambda)ds} dt
\end{eqnarray}
for each $n\ge 1,$ where the notation in (\ref{09June2020Eqn05}) is in use with $x_{(m)}$ being replaced by $y.$ This would yield the desired relation (\ref{09June2020Eqn07}) because by (\ref{Correction05})
\begin{eqnarray*}
&&f(T,x_{(m)},n-1)=\breve{P}_{(x,i)}^{\breve{\pi}}(x_{(m+1)}\in \Gamma,~\tau_{(m+1)}-\tau_{(m)}\in[0,T],\\
&&~ \mbox{exactly $n-1$ ficticious jumps over $[\tau_{(m)},\tau_{(m+1)}]$}|~x_{(0)},\tau_{(1)},x_{(1)},\dots,\tau_{(m)}-\tau_{(m-1)},x_{(m)}).
\end{eqnarray*}

Relation (\ref{09June2020Eqn06}) holds trivially when $n=1$ because
\begin{eqnarray*}
f(T,y,0)=\int_0^T \tilde{q}(\Gamma|\phi(y,t),\pi_m)e^{-\int_0^t (q_{\phi(y,s)}(\pi_m)+\lambda)ds}dt
\end{eqnarray*}
by definition.
Suppose (\ref{09June2020Eqn06}) holds. Then
\begin{eqnarray*}
&&f(T,y,n)=\int_0^T \lambda\int_0^{T-r_1}\int_{0}^{T-\sum_{i=1}^2r_i}\dots\int_{0}^{T-\sum_{i=1}^{n-1}r_i}\int_{0}^{T-\sum_{i=1}^{n}r_i} \tilde{q}(\Gamma|\phi(y,\sum_{i=1}^{n}r_i+t),\pi_m)\nonumber\\
&&\times e^{-\int_0^t (q_{\phi(y,\sum_{i=1}^{n}r_i+s)}(\pi_m)+\lambda)ds}\lambda^{n-1} e^{-\int_0^{\sum_{i=1}^{n}r_i} (q_{\phi(y,s)}(\pi_m)+\lambda)ds}   dt dr_{n}\dots dr_2d r_1\nonumber\\
&=&\int_0^T \lambda e^{-\int_0^{r_1} (q_{\phi(y,s)}(\pi_m)+\lambda)ds} \\
&&\times \left\{\int_0^{T-r_1}\int_{0}^{T-r_1-r_2}\dots \int_{0}^{T-r_1-\sum_{i=2}^{n}r_i} \tilde{q}(\Gamma|\phi(\phi(y,r_1),\sum_{i=2}^{n}r_i+t),\pi_m)\nonumber\right.\\
&&\left.\times e^{-\int_0^t (q_{\phi(\phi(y,r_1),\sum_{i=2}^{n}r_i+s)}(\pi_m)+\lambda)ds}\lambda^{n-1} e^{-\int_{0}^{\sum_{i=2}^{n}r_i} (q_{\phi(\phi(y,r_1),s)}(\pi_m)+\lambda)ds}   dt dr_{n}\dots dr_2 \right\}dr_1
\end{eqnarray*}
where the second equality holds because of (\ref{09June2020Eqn08}) and that
\begin{eqnarray*}
&&e^{-\int_0^{\sum_{i=1}^{n}r_i} (q_{\phi(y,s)}(\pi_m)+\lambda)ds} =e^{-\int_0^{r_1} (q_{\phi(y,s)}(\pi_m)+\lambda)ds} e^{-\int_{r_1}^{\sum_{i=1}^{n}r_i} (q_{\phi(y,s)}(\pi_m)+\lambda)ds}\\
&=&e^{-\int_0^{r_1} (q_{\phi(y,s)}(\pi_m)+\lambda)ds}
e^{-\int_{0}^{\sum_{i=2}^{n}r_i} (q_{\phi(\phi(y,r_1),s)}(\pi_m)+\lambda)ds}.
\end{eqnarray*}

Now after applying the induction supposition (\ref{09June2020Eqn06}) (with $y$ being replaced by $\phi(y,r_1)$, $T$ being replaced by $T-r_1$, and $\sum_{i=2}^n r_i$ playing the same role as $\sum_{i=1}^{n-1}r_i$ therein) to the previous inner integral, we see that
\begin{eqnarray*}
&&f(T,y,n)=\int_0^T \lambda e^{-\int_0^{r_1} (q_{\phi(y,s)}(\pi_m)+\lambda)ds }\\
&&\times \left\{\int_{0}^{T-r_1} \frac{\lambda^{n-1}t^{n-1}}{(n-1)!}\tilde{q}(\Gamma|\phi(\phi(y,r_1),t),\pi_m)e^{-\int_0^t (q_{\phi(\phi(x,r_1),s)}(\pi_m)+\lambda)ds}dt\right\}dr_1\\
&=&\int_0^T   \left\{\int_{0}^{T-r_1} \frac{\lambda^{n}t^{n-1}}{(n-1)!}\tilde{q}(\Gamma|\phi(y,r_1+t),\pi_m)e^{-\int_0^{t+r_1} (q_{\phi(x,s)}(\pi_m)+\lambda)ds}dt\right\}dr_1\\
&=&\int_0^T \int_0^t \frac{\lambda^n (t-r_1)^{n-1}}{(n-1)!} \tilde{q}(\Gamma|\phi(y,t),\pi_m)e^{-\int_0^{t} (q_{\phi(x,s)}(\pi_m)+\lambda)ds}dr_1 dt\\
&=&\int_0^T   \frac{\lambda^n t^{n}}{n!} \tilde{q}(\Gamma|\phi(y,t),\pi_m)e^{-\int_0^{t} (q_{\phi(x,s)}(\pi_m)+\lambda)ds}  dt,
\end{eqnarray*}
as desired, where the second equality follows from (\ref{09June2020Eqn08}) and that
\begin{eqnarray*}
e^{-\int_0^{r_1} (q_{\phi(y,s)}(\pi_m)+\lambda)ds } e^{-\int_0^t (q_{\phi(\phi(x,r_1),s)}(\pi_m)+\lambda)ds}=e^{-\int_0^{t+r_1} (q_{\phi(y,s)}(\pi_m)+\lambda)ds },
\end{eqnarray*}
and the third equality is by a change of variable and interchanging legitimately the order of integration. Part (a) is thus proved.

(b)
It follows from part (a) that \begin{eqnarray*}
&& \breve{V}((x,i_0),\breve{\pi})=\breve{E}^{\breve{\pi}}_{(x,i_0)}[e^{\int_0^\infty \breve{c}(\breve{\xi}_t,a)\breve{\pi}(da|\breve{\omega},t)dt}]\\
 &=&\breve{E}^{\breve{\pi}}_{(x,i_0)}\left[e^{\sum_{m=0}^\infty\int_{\tau_{(m)}}^{\tau_{(m+1)}}\int_A c(\phi(x_{(m)},s),a)\pi_{m}(da|x,\tau_{(1)},x_{(2)},\tau_{(2)}-\tau_{(1)},\dots,\tau_{(m)}-\tau_{(m-1)},x_{(m)},s-\tau_{(m)})ds} \right]\\
&=& {E}^{{\pi}}_{(x,i_0)}\left[e^{\sum_{m=0}^\infty\int_{t_{m}}^{t_{m+1}}\int_A c(\phi(x_{m},s),a)\pi_{m}(da|x,t_1,x_2,t_2-t_1,\dots,t_{m}-t_{m-1},x_{m},s-t_m)ds} \right]=V(x,\pi),
 \end{eqnarray*}
 where the second equality holds by (\ref{Correction05}) and the definition of $\breve{c}.$ The proof is completed. $\hfill\Box$

\begin{remark}\label{09June2020Remark01}
\begin{itemize}
\item[(a)] By Theorem \ref{09June2020Thm09}(b),
$\breve{V}^\ast(x)\le V^\ast(x)$ for each $x\in S.$
\item[(b)]
By inspecting the proof of Theorem \ref{09June2020Thm09} (see especially (\ref{09June2020Eqn10}) and (\ref{Correction05}) therein), one can tell that for a deterministic stationary policy in the auxiliary model, which depends on $(x,i)\in S\times \{-1,1\}$ only through $x\in S$, and is identified by a measurable mapping $\varphi$ from $S$ to $A$, $\breve{V}((x,i),\varphi)=V(x,\varphi)$ for all $x\in S$ and $i\in \{-1,1\}$. Therefore, if such a deterministic stationary policy $\varphi$ is optimal in the auxiliary model, then so is it in the original model, and $ V^\ast(x)=\breve{V}^\ast(x)=\breve{V}((x,i),\varphi)=V(x,\varphi)$ for each $x\in S.$
\item[(c)] By Theorem \ref{09June2020Thm09}(a) and the second equality in (\ref{09June2020Eqn15}), we see  $\overline{V}(x,\pi)=\breve{\overline{V}}((x,i),\breve{\pi})$ for the same policies $\pi$ and $\breve{\pi}$ in Theorem \ref{09June2020Thm09}, too.
\end{itemize}
\end{remark}

Let us introduce a natural set of compactness-continuity conditions on the original PDMDP.
\begin{condition}\label{GZyExponentialCondition02}
\begin{itemize}
\item[(a)] For each bounded measurable function $f$ on $S$ and each $x\in S$, $\int_S f(y)\tilde{q}(dy|x,a)$ is continuous in $a\in A.$
\item[(b)] For each $x\in S,$ the (nonnegative) function $c(x,a)$ is lower semicontinuous in $a\in A.$
\item[(c)] The action space $A$ is a compact Borel space.
\end{itemize}
\end{condition}

The usefulness of the auxiliary PDMDP also partially lies in the next observation.
\begin{lemma}\label{09June2020Lemma01}
If the original PDMDP model satisfies Conditions \ref{GZyExponentialCondition02} and \ref{GGZyExponentialConditionExtra}, then the auxiliary model satisfies the corresponding versions of Conditions \ref{GZyExponentialCondition02} and \ref{GGZyExponentialConditionExtra}, too.
\end{lemma}
\par\noindent\textit{Proof.} We only verify the version of Condition \ref{GZyExponentialCondition02}(a). For any bounded measurable function $f$ on $\breve{S}$,
it holds that
\begin{eqnarray*}
&&\int_{\breve{S}}f(y,j)\tilde{\breve{q}}(d(y,j)|(x,i),a)= \int_{\breve{S}}f(y,j){\breve{q}}(d(y,j)|(x,i),a)+f(x,i)\breve{q}_{(x,i)}(a)\\
&=& \int_S f(y,i)(q(dy|x,a)-\lambda\delta_x(dy))+\int_S f(y,-i)\lambda \delta_x(dy)+f(x,i)(\lambda+q_x(a))\\
&=&\int_S f(y,i)\tilde{q}(dy|x,a)-q_x(a)f(x,i)-\lambda f(x,i)+\lambda f(x,-i)+f(x,i)(\lambda+q_x(a))\\
&=&\int_S f(y,i)\tilde{q}(dy|x,a) +\lambda f(x,-i),
\end{eqnarray*}
which is clearly continuous in $a\in A$ when the original model satisfies Condition \ref{GZyExponentialCondition02}. $\hfill\Box$
\bigskip

The following statement was obtained in Theorem 3.1 and Remark 3.1 of \cite{GuoZhangAMO}.
\begin{proposition}\label{CorrectionRemark01}
Suppose Conditions \ref{GGZyExponentialConditionExtra} and \ref{GZyExponentialCondition02} are satisfied. In addition\footnote{See the footnote in Proposition \ref{09June2020PropositionA01}.}, $\inf_{(x,a)\in S\times A}q_x(a)>0$. Then the following assertions hold.
\begin{itemize}
\item[(a)] The value function $V^\ast$ is the minimal $[1,\infty]$-valued measurable solution to the following optimality equation:
  \begin{eqnarray}\label{09June2020Eqn11}
&& -(V(\phi(x,t))-V(x)) \\
&=&\int_0^t \inf_{a\in A}\left\{ \int_S V(y)\tilde{q}(dy|\phi(x,\tau),a)- (q_{\phi(x,\tau)}(a)-c(\phi(x,\tau),a) )V(\phi(x,\tau))\right\}d\tau \nonumber\\
&&\forall~t\in[0,\infty),x\in S^\ast:=\{x\in S:~V^\ast(x)<\infty\};~V(x)<\infty~\forall~x\in S^\ast; ~V(x)=\infty ~\forall~x\notin S^\ast.\nonumber
 \end{eqnarray}
In particular, $V^\ast(\phi(x,t))$ is absolutely continuous in $t$ for each $x\in S^\ast.$
 \item[(b)]  Any measurable mapping $\varphi$ from $S$ to $A$ such that
 \begin{eqnarray*}
 &&\inf_{a\in A}\left\{ \int_S V^\ast(y)\tilde{q}(dy|x,a)- (q_{x}(a)-c(x,a))V^\ast(x)\right\}\\
 &=&\int_S V^\ast(y)\tilde{q}(dy|x,\varphi(x))- (q_{x}(\varphi(x))-c(x,\varphi(x)))V^\ast(x),~\forall~x\in S^\ast.
 \end{eqnarray*}
 defines a deterministic stationary optimal policy in the original model. Such measurable selectors $\varphi$ exist.
\end{itemize}
\end{proposition}

As an application of Theorem \ref{09June2020Thm09} (more precisely, Remark \ref{09June2020Remark01} drawn from it), we may remove the redundant condition on the strong positivity of the transition intensities from Proposition \ref{CorrectionRemark01}.
\begin{corollary}\label{09June2020Corollary01}
Under Conditions \ref{GGZyExponentialConditionExtra} and \ref{GZyExponentialCondition02}, the assertions stated in Proposition \ref{CorrectionRemark01} all hold without requiring $\inf_{(x,a)\in S\times A}q_x(a)>0.$
\end{corollary}

\par\noindent\textit{Proof.} The statement follows from Remark \ref{09June2020Remark01} and applying Proposition \ref{CorrectionRemark01} to the auxiliary model, which is legitimate in view of Lemma \ref{09June2020Lemma01} and that the transition intensities in the auxiliary model are strongly positive. The details are as follows.

Step 1. We show that the value function $\breve{V}^\ast((x,i))$ in the auxiliary PDMDP model depends on $(x,i)$ only through $x$, and can thus be identified as $\breve{V}^\ast(x)$.

For this, we will apply the following result from \cite{GuoZhangAMO}: under the conditions in Proposition \ref{CorrectionRemark01}, including that the transition intensities are strongly positive:
\begin{itemize}
\item The value function $V^\ast$ in the original model is the minimal $[1,\infty]$-valued measurable solution to the optimality equation $V={\cal T}\circ V$, where
\begin{eqnarray}\label{Correction01}
{\cal T}\circ V(x)&:=& \inf_{\rho\in {\cal R}}\left\{\int_0^\infty e^{-\int_0^\tau (q_{\phi(x,s)}(\rho_s)-c(\phi(x,s),\rho_s))ds} \left(\int_S V(y)\tilde{q}(dy|\phi(x,\tau),\rho_\tau)\right)d\tau\right.\nonumber\\
&&\left. +e^{-\int_0^\infty q_{\phi(x,s)}(\rho_s)ds}e^{\int_0^\infty c(\phi(x,s),\rho_s)ds} \right\},~\forall~x\in S,
\end{eqnarray}
Here and below, ${\cal R}$ is the space of $\mathbb{P}(A)$-valued measurable mappings $\rho=(\rho_t(da))$ on $(0,\infty)$\footnote{Two elements in ${\cal R}$ that coincide almost everywhere are not distinguished.},  and $e^{-\int_0^\infty q_{\phi(x,s)}(\rho_s)ds}e^{\int_0^\infty c(\phi(x,s),\rho_s)ds}:=0$
whenever
$e^{-\int_0^\infty q_{\phi(x,s)}(\rho_s)ds}=0.$
\item The value function $V^\ast$ can be obtained from the successive approximation: $V^\ast(x)=\lim_{n\rightarrow \infty}{\cal T}\circ V_0(x)$ with $V_0(x)\equiv 1$.
\end{itemize}

According to Lemma \ref{09June2020Lemma01} and that the transition intensities in the auxiliary model are strongly positive (as $\inf_{x\in S,a\in A}\breve{q}_{(x,i)}(a)\ge \lambda>0$), we may apply the result just quoted above to the auxiliary model and conclude that $\breve{V}^\ast$ is the minimal $[1,\infty]$-valued measurable function to the following equation
\begin{eqnarray*}
&&\breve{V}((x,i))\\
&=&\inf_{\rho\in{\cal R}}\left\{\int_0^\infty e^{\int_0^\theta \breve{c}(\breve{\phi}((x,i),s),\rho_s)ds} e^{-\int_0^\theta \breve{q}_{\breve{\phi}((x,i),s)}(\rho_s) ds}\left(\int_{\breve{S}} \breve{V}((y,j))\tilde{\breve{q}}(d(y,j)|\breve{\phi}((x,i),\theta),\rho_\theta) \right)d\theta\right\}\\
&=&\inf_{\rho\in{\cal R}}\left\{\int_0^\infty e^{\int_0^\theta c(\phi(x,s),\rho_s)ds} e^{-\int_0^\theta (q_{\phi(x,s)}(\rho_s)+\lambda )ds}\right.\\
&&\left.\times\left(\int_S \tilde{q}(dy|\phi(x,\theta),\rho_\theta)\breve{V}((y,i))+\lambda \breve{V}((\phi(x,\theta),-i))\right)d\theta\right\}.
\end{eqnarray*}
Moreover, $\breve{V}^\ast$ is the pointwise limit of the sequence of functions $\{\breve{V}_n\}_{n=0}^\infty$ with
\begin{eqnarray*}
&&\breve{V}_0((x,i)):\equiv 1,\\
&&\breve{V}_{n+1}((x,i)):=\inf_{\rho\in{\cal R}}\left\{\int_0^\infty e^{\int_0^\theta c(\phi(x,s),\rho_s)ds} e^{-\int_0^\theta (q_{\phi(x,s)(\rho_s)}+\lambda )ds}\right.\\
&&\left.\left(\int_S \tilde{q}(dy|\phi(x,\theta),\rho_\theta)\breve{V}_n((y,i))+\lambda \breve{V}_n((\phi(x,\theta),-i))\right)d\theta\right\}.
\end{eqnarray*}
An inductive argument reveals that $\breve{V}_{n+1}((x,i))$ does not depend on $i$ for all $n\ge 0$ and thus $\breve{V}^\ast((x,i))$ does not depend on $i.$  Below, we write $\breve{V}^\ast(x)$ for $\breve{V}^\ast((x,i))$.

Step 2. Again by Lemma \ref{09June2020Lemma01} and that the transition intensities in the auxiliary model are strongly positive,  we apply Proposition \ref{CorrectionRemark01}(b) to the auxiliary model to obtain a deterministic stationary optimal policy $\varphi$. It is possible to take $\varphi$, which only depends on $x\in S$ (independent on $i\in \{-1,1\}$) because for each $(x,i)\in \breve{S}^\ast:=\{x\in S: \breve{V}^\ast(x)<\infty\}\times \{-1,1\}$,
\begin{eqnarray*}
 &&\int_{\breve{S}} \breve{V}^\ast(y)\tilde{\breve{q}}(d(y,j)|(x,i),a)- (\breve{q}_{(x,i)}(a)-\breve{c}((x,i),a))\breve{V}^\ast(x))\\
 &=&  \int_S \breve{V}^\ast(y)\tilde{q}(dy|x,a)- (q_{x}(a)-c(x,a))\breve{V}^\ast(x))
 \end{eqnarray*}
does not involve $i\in\{-1,1\}$, where the equality holds by the definition of $\breve{q}$ and $\breve{c}$ and a similar calculation as the one in the proof of Lemma \ref{09June2020Lemma01}.

Step 3. Step 2 and Remark \ref{09June2020Remark01}(b) imply $\breve{V}^\ast(x)=V^\ast(x)$. Note that the optimality equations for both the original model and the auxiliary model are the same and given by (\ref{09June2020Eqn11}), the statement of this corollary follows from applying again Proposition \ref{CorrectionRemark01} to the auxiliary model. $\hfill\Box$

\appendix
\section{Appendix: Induced DTMDP}\label{NoteGZSec02}
We shall formulate a DTMDP model induced from the PDMDP model with total cost criteria by inspecting the PDMDP at each of its jump moments and regarding the relaxed control functions used during a sojourn time as the actions in the DTMDP. The first coordinate in the state space of the induced DTMDP records the most recent sojourn time, and the second coordinate records the state in the PDMDP immediately after the corresponding jump. This is to serve the formulation of Example \ref{JapanExample01} below.

The DTMDP induced by the PDMDP $\{S,A,q,\phi,c\}$ is specified by the following system primitives:
\begin{itemize}
\item The state space is $\textbf{X}:=((0,\infty)\times S)\bigcup\{(\infty,x_\infty)\}$. Whenever the topology is concerned, $(\infty,x_\infty)$ is regarded as an isolated point in $\textbf{X}.$
\item The action space is $\textbf{A}:={\cal R}$, where ${\cal R}$ was defined in the proof of Corollary \ref{09June2020Corollary01}. We endow ${\cal R}$ with the Young topology.\footnote{The Young topology on ${\cal R}$ is the weakest
topology with respect to which the function
$
\rho\in {\cal{R}}\rightarrow \int_0^\infty \int_A
f(t,a)\rho_t(da)dt
$
is continuous for each strongly integrable Carath\'{e}odory function
$f$ on $(0,\infty)\times A$ . Here a real-valued measurable
function $f$ on $(0,\infty)\times A$ is called a strongly
integrable Carath\'{e}odory function if for each fixed
$t\in(0,\infty)$, $f(t,a)$ is continuous in $a\in A,$ and for each
fixed $a\in A,$ $\sup_{a\in A}|f(t,a)|$ is integrable in $t$,
i.e., $\int_0^\infty \sup_{a\in A}|f(t,a)|dt<\infty.$}
\item The transition kernel $p$ on ${\cal B}(\textbf{X})$ from $\textbf{X}\times \textbf{A}$  is given for each $\rho=(\rho_t(da))_{t>0}\in \textbf{A}$ by
    \begin{eqnarray*}
    p(\Gamma_1\times\Gamma_2|(\theta,x),\rho)&:=&\int_{\Gamma_2} e^{-\int_0^t q_{\phi(x,s)}(\rho_s)ds}\tilde{q}(\Gamma_1|\phi(x,t),\rho_t)dt,\nonumber\\
    &&~\forall~\Gamma_1\in {\cal B}(S),~\Gamma_2\in {\cal B}((0,\infty)),~x\in S,~\theta\in (0,\infty),\nonumber\\
     p(\{(\infty,x_\infty)\}|(\theta,x),\rho)&:=&e^{-\int_0^\infty q_{\phi(x,s)}(\rho_s)ds},~\forall~x\in S,~\theta\in (0,\infty);\nonumber\\
     p(\{(\infty,x_\infty)\}|(\infty,x_\infty),\rho)&:=&1.\nonumber
\end{eqnarray*}
(Recall that the notation $q(dy|x,\rho_t)=\int_A q(dy|x,a)\rho_t(da)$ is in use.)
\item The cost function $l$ is a $[0,\infty]$-valued measurable function on $\textbf{X}\times\textbf{A}\times \textbf{X}$ given by
\begin{eqnarray*}\label{ZyExponential55}
l((\theta,x),\rho,(\tau,y)):=\int_0^\infty I\{s<\tau\} c(\phi(x,s),\rho_s)ds,~\forall~((\theta,x),\rho,(\tau,y))\in \textbf{X}\times\textbf{A}\times \textbf{X}.
\end{eqnarray*}
\end{itemize}

For the induced DTMDP $\{\textbf{X},\textbf{A},p,l\}$, following the reasoning in the proof of Lemma 3.2 of \cite{Costa:2013} and Chapter 4 of \cite{Davis:1993}, one can see the following statement.
\begin{proposition}\label{09June2020PropositionA01}
Under Conditions \ref{GGZyExponentialConditionExtra} and \ref{GZyExponentialCondition02}, for each $(\theta,x)\in \textbf{X}$ and $(\tau,y)\in\textbf{X}$, $\rho\in\textbf{A}\rightarrow l((\theta,x),\rho,(\tau,y))$ is lower semicontinuous, and $\textbf{A}$ is a compact Borel space. If in addition, the transition intensities are strongly positive\footnote{This requirement was unfortunately missing and overlooked in \cite{GuoZhangAMO}. Indeed, the proof of Lemma 4.1 of \cite{GuoZhangAMO} made use of the strong Feller property of the transition kernel $p$ in the induced DTMDP, which could fail to hold without this additional requirement, as demonstrated in Example \ref{JapanExample01}. Since the rest of the arguments in \cite{GuoZhangAMO} are largely based on that lemma, this missing requirement on the strong positivity of transition intensities should be added in \cite{GuoZhangAMO} wherever appropriate. We thank Dr Yonghui Huang (Sun Yat-Sen University, China) for drawing our attention on this inaccuracy.}, then for each $(\theta,x)\in \textbf{X}$, the function $\rho\in \textbf{A}\rightarrow \int_{\textbf{X}}f(z)p(dz|(\theta,x),\rho)$ is continuous for each bounded measurable function $f$ on $\textbf{X}$.
\end{proposition}

The next example shows that if the transition intensities are not strongly positive, then it can happen that $\rho\in \textbf{A}\rightarrow \int_{\textbf{X}}f(z)p(dz|(\theta,x),\rho)$ is not continuous for  some bounded measurable function $f$ on $\textbf{X}$.

\begin{example}\label{JapanExample01}
Suppose $S$ is any finite set (endowed with discrete topology), and $A=[0,1]$, which is a compact Borel space, $q_x(a)=a$ and $c(x,a)\equiv 0$, and $\phi(x,t)\equiv x$. Evidently, Conditions \ref{GGZyExponentialConditionExtra} and \ref{GZyExponentialCondition02} are satisfied by this PDMDP model. Consider $\rho\in \textbf{A}$ and $(\rho^{(n)})\subseteq \textbf{A}$ defined as follows: for each $t\ge 0,$ $\rho^{(n)}_t(da)=\delta_{\frac{1}{n}}(da)$, and $\rho_t(da)=\delta_0(da).$
Then for each strongly integrable Carath\'{e}odory function $g(t,a)$,
\begin{eqnarray*}
\int_0^\infty g(t,\rho^{(n)}_t)dt-\int_0^\infty g(t,\rho^{(0)}_t)dt=\int_0^\infty (g(t,\frac{1}{n})-g(t,0))dt\rightarrow 0
\end{eqnarray*}
as $n\rightarrow \infty,$ by using the dominated convergence theorem. Thus, $\rho^{(n)}\rightarrow \rho$ as $n\rightarrow \infty.$  (Recall that $\textbf{A}$ is endowed with the Young topology.)
Now for $f(t,x)\equiv 0$ on $(0,\infty)\times S$ and $f(\infty,x_\infty)=1$,
\begin{eqnarray*}
  \int_{\textbf{X}}f(z)p(dz|(\theta,x),\rho^{(n)})=e^{-\int_0^\infty q_x(\rho^{(n)}_s)ds}=e^{-\int_0^\infty \frac{1}{n}ds}=0<  1=e^{-\int_0^\infty 0 ds}= \int_{\textbf{X}}f(z)p(dz|(\theta,x),\rho).
\end{eqnarray*}
\end{example}


\end{document}